\documentclass[12pt]{amsart}
\usepackage[latin1]{inputenc}
\usepackage{amscd}
\usepackage{latexsym}
\usepackage{pstcol,pst-plot,pst-3d}
\usepackage{mathptmx}
\usepackage{multicol}

\psset{unit=0.7cm,linewidth=0.8pt,arrowsize=2.5pt 4}

\newpsstyle{fatline}{linewidth=1.5pt}
\newpsstyle{fyp}{fillstyle=solid,fillcolor=verylight}
\definecolor{verylight}{gray}{0.97}
\definecolor{light}{gray}{0.9}
\definecolor{medium}{gray}{0.85}

\def\NZQ{\Bbb}               

\def\PP{{\NZQ P}}
\def\FFF{{\NZQ F}}
\def\TT{{\NZQ T}}

\def\opn#1#2{\def#1{\operatorname{#2}}} 
\opn\chara{char} \opn\length{\ell} \opn\pd{pd} \opn\rk{rk}
\opn\projdim{proj\,dim} \opn\injdim{inj\,dim} \opn\rank{rank}
\opn\depth{depth} \opn\grade{grade} \opn\height{height}
\opn\embdim{emb\,dim} \opn\codim{codim}

\opn\Tr{Tr} \opn\bigrank{big\,rank}
\opn\superheight{superheight}\opn\lcm{lcm}
\opn\trdeg{tr\,deg}%
\opn\reg{reg} \opn\lreg{lreg} \opn\skel{skel}
\opn\multideg{multideg}
\opn\div{div} \opn\Div{Div} \opn\cl{cl} \opn\Cl{Cl}
\opn\Spec{Spec} \opn\Supp{Supp} \opn\supp{supp} \opn\Sing{Sing}
\opn\Ass{Ass}
\opn\Ann{Ann} \opn\Rad{Rad} \opn\Soc{Soc}
\opn\Ker{Ker} \opn\Coker{Coker} \opn\Im{Im} \opn\Hom{Hom}
\opn\Tor{Tor} \opn\Ext{Ext} \opn\End{End} \opn\Aut{Aut}
\opn\id{id}

\opn\nat{nat}
\opn\pff{pf}
\opn\Pf{Pf} \opn\GL{GL} \opn\SL{SL} \opn\mod{mod} \opn\ord{ord}
\opn\aff{aff} \opn\con{conv} \opn\relint{relint} \opn\st{st}
\opn\lk{lk} \opn\cn{cn} \opn\core{core} \opn\vol{vol}
\opn\link{link} \opn\star{star} \opn\skel{skel}
\opn\gr{gr}

\def\pot#1#2{#1[\kern-0.28ex[#2]\kern-0.28ex]}

\opn\dirlim{\underrightarrow{\lim}}
\opn\inivlim{\underleftarrow{\lim}}
\let\union=\cup
\let\sect=\cap
\let\dirsum=\oplus

\let\iso=\cong
\let\Union=\bigcup
\let\Sect=\bigcap

\let\to=\rightarrow
\let\To=\longrightarrow
\def\Implies{\ifmmode\Longrightarrow \else
     \unskip${}\Longrightarrow{}$\ignorespaces\fi}
\def\implies{\ifmmode\Rightarrow \else
     \unskip${}\Rightarrow{}$\ignorespaces\fi}
\def\iff{\ifmmode\Longleftrightarrow \else
     \unskip${}\Longleftrightarrow{}$\ignorespaces\fi}

\let\:=\colon
\newtheorem{Theorem}{Theorem}[section]
\newtheorem{Lemma}[Theorem]{Lemma}
\newtheorem{Corollary}[Theorem]{Corollary}
\newtheorem{Proposition}[Theorem]{Proposition}
\newtheorem{Remark}[Theorem]{Remark}

\newtheorem{Definition}[Theorem]{Definition}

\let\epsilon\varepsilon
\let\phi=\varphi
\let\kappa=\varkappa
\def\qed{\ifhmode\textqed\fi
   \ifmmode\ifinner\quad\qedsymbol\else\dispqed\fi\fi}
\def\textqed{\unskip\nobreak\penalty50
    \hskip2em\hbox{}\nobreak\hfil\qedsymbol
    \parfillskip=0pt \finalhyphendemerits=0}
\def\dispqed{\rlap{\qquad\qedsymbol}}

\def\BB{{\mathcal B}}

\def\LL{{\mathcal L}}
\def\II{{\mathcal I}}
\def\JJ{{\mathcal J}}
\def\SS{{\mathcal S}}

\opn\inii{in} \opn\inim{inm} \opn\rate{rate}
\begin{document}
\title{Monomial ideals arising from distributive lattices}
\author{Xinxian Zheng}
\address{Xinxian Zheng, Fachbereich Mathematik und
Informatik, Universit\"at Duisburg-Essen, 45117 Essen, Germany}
\email{xinxian.zheng@uni-essen.de}
\date{}
\maketitle

\begin{abstract}
The free resolution and the Alexander dual of squarefree monomial
ideals associated with certain subsets of  distributive lattices
are studied.
\end{abstract}

\section*{Introduction}
Let $\LL$ be a finite distributive lattice. By Birkhoff's
fundamental structure theorem, there is a unique poset (partially
ordered set) $P$ such that $\LL$ is isomorphic to the poset
$\JJ(P)$ consisting of all poset ideals (including the empty set)
of $P$, ordered by inclusion. In fact, $P$ can be chosen as the
set of all join-irreducible elements of $\LL$. Let $K$ be a field
and $S = K[\{ x_p, y_p \}_{p \in P}]$ the polynomial ring in
$2|P|$ variables over $K$ with $\deg x_p = 1$ and $\deg y_p=1$ for
all $p\in P$, and let $\SS\subset \LL$ be any subset of $\LL$. The
Hibi ideal $H_\SS$ associated with $\SS$ is the monomial ideal in
$S$ generated by the monomials $u_p$ with $p\in\SS$, where
$u_p=x_{\ell(p)}y_{P\setminus \ell(p)}$ and where $\ell(p)$ is the
principal poset ideal $\{q\in P\: q\leq p\}$ in $P$.

In \cite{HHZ2} it is shown that for any poset ideal $\II$ of
$\LL$, the Hibi ideal $H_\II$ has a linear resolution. In this
article, we consider more generally the ideal $H_\SS$ where $\SS$
is a segment of $\LL$ (see Definition \ref{segment}). For example,
any poset ideal $\II$, or any poset coideal $\JJ$ of $\LL$, as
well as their intersection are segments in $\LL$. In the third
section we describe in Theorem \ref{equal} when $H_\II\sect
H_\JJ=H_{\II\sect\JJ}$, and in Theorem \ref{linear} it is said
when this ideal has a linear resolution. In particular this
answers a question which was raised in \cite{HHZ2}, see Corollary
\ref{one line} and \ref{planar}. We also show in Theorem
\ref{empty} that the ideal $H_\II\sect H_\JJ$ has always a linear
resolution, if $\II\union\JJ=\LL$ and $\II\sect\JJ=\emptyset$.

Let $G$ be a Cohen--Macaulay bipartite graph on the vertex set
$V\union V'$ with $V\sect V'=\emptyset$ and $|V|=|V'|=n$, and
$S=K[x_1,\ldots,x_n,y_1,\ldots,y_n]$ the polynomial ring over a
field $K$. In \cite[Theorem 2.4]{HH} the authors showed that the
vertices $V=\{x_1,\ldots, x_n\}$ and $V'=\{y_1,\ldots, y_n\}$ can
be labeled such that there exists a partial order $<$ on $V$ with
the property that $\{x_i,y_j\}$ is an edge of $G$ if and only if
$x_i\leq x_j$. Moreover it is shown that for $P = (V, <)$ the
distributive lattice $\JJ(P)$ satisfies $H_{\JJ(P)}^*=I(G)$. Here,
for any subset $\SS$ of $\JJ(P)$ we denote by  $H^*_\SS$  the
defining ideal of the Stanley--Reisner of the Alexander dual of
$\Gamma$, where $\Gamma$ is the simplicial complex defined by the
equation $H_\SS=I_\Gamma$.

Later, in \cite{HHZ2} the authors considered more generally
simplicial complexes $\Delta$ on the vertex set $V\union V'$ with
$V\sect V'=\emptyset$ and $|V|=|V'|$ such that
\begin{itemize}
\item[(1)] there is no $F\in \mathcal{F}(\Delta)$ with $F\subset
V$,

\item[(2)] $G=\{F\in\mathcal{F}(\Delta)\: F\sect V\neq
\emptyset,\quad F\sect V'\neq \emptyset\}$ is a Cohen--Macaulay
bipartite graph with no isolated vertex,
\end{itemize}
and  studied when the facet ideal $I(\Delta)$ of $\Delta$ is
Cohen-Macaulay. As a further  generalization we consider in the
second section of this article simplicial complexes $\Delta$
satisfying only condition (2), and show (Theorem \ref{unmixed})
that $\Delta$ is unmixed and each minimal vertex cover of $\Delta$
has cardinality $n$ if and only if there exists a segment $\SS$ of
some distributive lattice $\LL$ such that $H^*_\SS=I(\Delta)$.

I would like to thank J\"urgen Herzog for many helpful comments
and discussions.

\section{Preparations}
In this section we recall some basic facts on lattices and
simplicial complexes and fix some notation. As a general reference
for posets and lattices we refer the reader to \cite{H} and
\cite{S},  and to \cite{S2}, \cite{BH} and \cite{Z} concerning
simplicial complexes, Stanley-Reisner  and facet ideals.

Let $P$ be any finite poset (partially ordered set), and let
$\alpha,\beta\in P$ with $\alpha\leq \beta$. The set
\[
[\alpha,\beta]=\{\gamma\in P:\alpha\leq \gamma\leq \beta\}
\]
is called the {\em interval} between $\alpha$ and $\beta$ in $P$.

Let $P$ be a poset and $\alpha,\beta\in P$. If $\alpha<\beta$ and
for each element $\gamma\in P$ with $\alpha\leq\gamma\leq\beta$,
we have either $\gamma=\alpha$ or $\gamma=\beta$, then we say
$\beta$ {\em covers} $\alpha$, or $\alpha$ is a {\em lower
neighbor} of $\beta$, or  $\beta$ is an {\em upper neighbor} of
$\alpha$.

An element in a poset $P$ may have more than one upper neighbor
(resp.\ lower neighbor) or have no upper neighbor (resp.\ lower
neighbor). An element in a poset $P$ which has exactly one lower
neighbor is called a {\em join irreducible element} of $P$. The
set of all join irreducible elements with the induced order is a
poset, called the {\em join irreducible subposet} of $P$.
Conversely, an element in a poset $P$ which has exactly one upper
neighbor is called a {\em meet irreducible element} of $P$.

A {\em chain} is a poset in which any two elements are comparable.
A subset $C$ of a poset $P$ is called a {\em chain} if $C$ is a
chain when regarded as a subposet of $P$. The {\em length}
$\ell(C)$ of a finite chain is defined by $\ell(C)=|C|-1$. The
{\em length} (or {\em rank}) of a finite poset $P$ is
$\ell(P):=\max\{\ell(C):\, C \text{ is a chain of } P\}$. If every
maximal chain of $P$ has the same length $r$, then we say $P$ is
{\em graded of rank} $r$. In this case there is a unique {\em rank
function} $\rho: P\to \{0,\ldots,r\}$ such that $\rho(\alpha)=0$
if $\alpha$ is a minimal element of $P$, and
$\rho(\beta)=\rho(\alpha)+1$ if $\beta$ covers $\alpha$ in $P$. If
$\rho(\alpha)=i$, then we say $\alpha$ has {\em rank} $i$.

Later, we need the {\em dual poset} of  $P$. This is the poset
$\widetilde{P}$ on the same set as $P$, but such that $\alpha\leq
\beta$ in $\widetilde{P}$ if and only if $\beta\leq \alpha$ in
$P$.

A {\em lattice} is a poset ${\mathcal L}$ for which each pair of
elements $\alpha$ and $\beta$ has a least upper bound (called the
{\em join} of $\alpha$ and $\beta$, denoted by $\alpha\vee \beta$)
and a greatest lower bound (called the {\em meet} of $\alpha$ and
$\beta$, denoted by $\alpha\wedge\beta$).

One sees immediately from the definition that in a lattice
${\mathcal L}$, there is a unique element $\mu$ satisfies that
$\mu\geq \alpha$ for any $\alpha\in {\mathcal L}$.  This element
is called  the {\em maximum} of ${\mathcal L}$, and  denoted by
$\hat {1}$. Similarly, there is a unique element $\nu$ satisfies
$\nu\leq \alpha$ for any $\alpha\in {\mathcal L}$. This element is
called  the {\em minimum} of ${\mathcal L}$, and denoted by $\hat
{0}$.

A {\em poset ideal (coideal)} of a poset $P$ is a subset $I$ of
$P$ such that if $\alpha\in I$ and $\beta<\alpha$
($\beta>\alpha$), then $\beta\in I$. The maximal (minimal)
elements in $I$ are called the {\em generators} of $I$. The set of
generators is denoted by $G(I)$.

\begin{Remark}
\label{veryeasy} {\em Let $I\subset P$. Then the following
conditions are equivalent:
\begin{enumerate}
\item[(1)]  $I$ is a poset ideal (coideal) in $P$;

\item[(2)] $P\setminus I$ is a poset coideal (ideal) of $P$;

\item[(3)] $\widetilde{I}$ is a poset coideal (ideal) of
$\widetilde{P}$.
\end{enumerate}}
\end{Remark}

Let $P$ be an arbitrary finite poset and write ${\mathcal J}(P)$
for the poset which consists of all poset ideals of $P$ ordered by
inclusion.

For example, if $P$ is an antichain, i.e., any two elements of $P$
are incomparable, then ${\mathcal J}(P)\cong \BB_P$, where $\BB_P$
is the Boolean lattice consisting of all subsets of $P$. The rank
of $\BB_P$ is the cardinality of $P$.

Since the union $I\cup J$ and the intersection $I\cap J$ of poset
ideals $I$ and $J$ of $P$ are also poset ideals of $P$, the poset
${\mathcal J}(P)$ is in fact a lattice.

We say that a poset $P$ is {\em isomorphic} to a poset $Q$ if
there exists a bijection $\theta: P\rightarrow Q$ such that
$\alpha\leq \beta$ in $P$ if and only if
$\theta(\alpha)\leq\theta(\beta)$ in $Q$.

The most important class of lattices from the combinatorial point
of view is the distributive lattice. And one of the most
influential results in the classical lattice theory is Birkhoff's
fundamental structure theorem for the finite distributive lattice.

\begin{Theorem}[Birkhoff]
\label{Birkhoff} Let ${\mathcal L}$ be a finite distributive
lattice. Then there exists a unique (up to isomorphism) poset $P$
such that ${\mathcal L}$ is isomorphic to ${\mathcal J}(P)$.
\end{Theorem}

One finds the proof, for example, in \cite[Theorem 3.4.1]{S}. In
fact, $P$ can be chosen as the join irreducible subposet of
${\mathcal L}$.

The outline of the proof is as follows: Let $\LL$ be a finite
lattice, and let $P$ be the set of join irreducible elements of
$\LL$. As in \cite{HHZ2} we associate to each element $p\in \LL$
the poset ideal $\ell(p)=\{q\in P\: q\leq p\}$ of $P$. This
establishes a map $\ell\: \LL\to \JJ(P)$, which we call the {\em
canonical embedding} into the distributive lattice $\JJ(P)$. Note
that $\ell$ is an isomorphism if and only if $\LL$ is
distributive.

We call the cardinality of $\ell(p)$ the {\em degree} of $p$, and
denote it by $\deg p$.

The map $\ell$ has the following properties.

\begin{Lemma}
\label{better} Let ${\mathcal L}$ be a finite lattice, $\ell$ the
canonical embedding and  $s,t\in {\mathcal L}$ any two elements.
We have
\begin{enumerate}
\item[(1)]  $s=t$ if and only if $\ell(s)= \ell(t)$;

\item[(2)] $s\leq t$ if and only if $\ell(s)\subseteq \ell(t)$;

\item[(3)] $\ell(s)\cap \ell(t)=\ell(s\wedge t)$.
\end{enumerate}
\end{Lemma}

This lemma implies in particular that $\ell$ is an injective order
preserving map. In general however, $\ell$ is not an embedding of
lattices.

As a consequence of Remark \ref{veryeasy} we have the following
\begin{Lemma}
\label{duallattice} $\widetilde{\LL}\iso \JJ(\widetilde{P})$.
\end{Lemma}
\begin{proof}
Let $q\in \widetilde{\LL}$. Since the underlying set of
$\widetilde{\LL}$ is the same as that of $\LL$, we may apply
$\ell\: \LL\to \JJ(P)$ to $q$. Then $\widetilde{\LL}\to
\JJ(\widetilde{P})$, $q\mapsto P\setminus \ell(q)$ is the desired
isomorphism.
\end{proof}

We now introduce the squarefree monomial ideal $H_\LL$ associated
with a finite lattice $\LL$. Let $K$ be a field and
$S=K[\{x_p,y_p\}_{p\in P}]$  the polynomial ring in $2|P|$
variables over $K$. For each element $q\in\LL$ write
\[
u_q=\prod_{p\in \ell(q)}x_p\prod_{p\in P\setminus \ell(q)}y_p,
\]
and set $H_\LL=(u_q)_{q\in\LL}$. We call $H_\LL$ the {\em Hibi
ideal} of $\LL$. It is easy to see that the height of $H_\LL$ is
2.

\medskip
Recall that a finite lattice ${\mathcal L}$ is {\em  upper
semimodular} if ${\mathcal L}$ satisfies either of the following
two conditions.
\begin{enumerate}
\item[(1)] ${\mathcal L}$ is graded, and the rank function $\rho$
of ${\mathcal L}$ satisfies
$\rho(\alpha)+\rho(\beta)\geq\rho(\alpha\wedge
\beta)+\rho(\alpha\vee \beta)$ for all $\alpha,\beta\in{\mathcal
L}$.

\item[(2)] If $\alpha$ and $\beta$ both cover $\alpha\wedge
\beta$, then $\alpha\vee \beta$ covers both $\alpha$ and $\beta$.
\end{enumerate}

Among the upper semimodular lattices the distributive lattices can
be characterized algebraically. For this we need following
concept: let $I$ be a monomial ideal in a polynomial ring with the
(unique) minimal set $G(I)$ of monomial generators. The ideal $I$
is called a {\em linear quotient ideal} if the elements of $G(I)$
can be ordered $u_1, \ldots, u_m$ such that the colon ideals
$(u_1,\ldots, u_{i-1}):u_{i}$ are generated by variables. If $I$
is squarefree, then $I=(u_1, \ldots, u_m)$ is a linear quotient
ideal (in this order) if and only if for each $i$ and each $j<i$
there exists $k<i$ such that $u_k/[u_k,u_i]$ is a variable which
divides $u_j$. Here $[u,v]$ denotes the greatest common divisor of
$u$ and $v$.

It is easy to see that a linear quotient ideal $I$ has a linear
resolution, if all generators of $I$ have the same degree. If $I$
is an ideal with linear resolution, then all generators of $I$
have the same degree, say $d$. In this case we also say $I$ has a
$d$-linear resolution.

The following  characterization of finite distributive lattices is
an immediate consequence of \cite[Theorem 1.3]{HHZ2}.

\begin{Proposition}
\label{main} Let $\LL$ be an arbitrary finite upper semimodular
lattice. The following conditions are equivalent:
\begin{enumerate}
\item[(1)] $\LL$ is distributive;

\item[(2)] $H_\LL$ has linear quotients;

\item[(3)] $H_\LL$ has a linear resolution;

\item[(4)]  $H_\LL$ has linear relations.
\end{enumerate}
\end{Proposition}

\medskip Now we recall some concepts related to simplicial complex and
fix some notation.
 Let $\Delta$ be a simplicial complex on the vertex set
$[n]=\{1,\ldots,n\}$, $R=K[x_1,\ldots,x_n]$ the polynomial ring in
$n$ variables over a field $K$. We denote ${\mathcal F}(\Delta)$
the set of facets (maximal faces) of $\Delta$.  The simplicial
complex
\[
\Delta^\vee=\{[n]\setminus F\: F\not\in\Delta\}
\]
is called the {\em Alexander dual}  of $\Delta$. One has
$(\Delta^\vee)^\vee=\Delta$.

A {\em vertex cover} of $\Delta$ is a set $G\subset [n]$ such that
$G\sect F \neq\emptyset$ for all $F\in {\mathcal F}(\Delta)$. We
say a vertex cover $G$ of $\Delta$ is {\em minimal}, if each
proper subset of $G$ is not a vertex cover of $\Delta$. We denote
by ${\mathcal C}(\Delta)$ the set of minimal vertex covers of
$\Delta$. If all the minimal vertex cover of $\Delta$ have the
same cardinality, then we say $\Delta$ is {\em unmixed}. For
$F=\{i_1,\ldots, i_k\}\subset [n]$, let $P_F$ be the prime ideal
generated by $x_{i_1},\ldots, x_{i_k}$, and set $F^c=[n]\setminus
F$. As usual we denote by $I_\Delta$ the Stanley--Reisner ideal
and $K[\Delta]=R/I_\Delta$ the Stanley--Reisner ring of $\Delta$.
The following proposition is an elementary but important property
of the Stanley--Reisner ideal. One finds the proof for example in
\cite[Theorem 5.1.4]{BH}.

\begin{Proposition}
\label{prime} Let $\Delta$ be a simplicial complex over the vertex
set $[n]$. Then $$I_\Delta=\Sect_{F\in{\mathcal F}(\Delta)}
P_{F^c}.$$
\end{Proposition}

 The {\em facet ideal} is defined to be
\[
I(\Delta)=(x_F\: F\in{\mathcal F}(\Delta)),
\]
where $x_F=\prod_{i\in F}x_i$.

Let $\Gamma$ be the unique simplicial complex such that
$I_\Delta=I(\Gamma)$. Then
\begin{eqnarray}
\label{cover} I_{\Delta}=\Sect_{F\in {\mathcal C}(\Gamma)}
P_F\quad \text{and} \quad I_{\Delta^\vee} = (x_F\: F\in {\mathcal
C}(\Gamma)).
\end{eqnarray}

Set $\Delta^c=\langle F^c\: F\in {\mathcal F}(\Delta)\rangle.$
Then
\begin{eqnarray}
\label{formula1}
 I_{\Delta^\vee}=I(\Delta^c).
\end{eqnarray}
The easy proofs can be found  for example in \cite{HHZ1}.

The proof of the following  simple lemma can be found for example
in \cite[Proposition 1.8]{F}.

\begin{Lemma}
\label{minimalcover} Let $\Delta$ be a simplicial complex on the
vertex set $[n]$ and $I(\Delta)$ the facet ideal of $\Delta$. Then
an ideal $P=(x_{i_1},\ldots,x_{i_s})$ is a minimal prime of
$I(\Delta)$ if and only if $\{i_1,\ldots,i_s\}$ is a minimal
vertex cover of $\Delta$.
\end{Lemma}

Let $\Delta$ be a simplicial complex and $\Gamma$ the unique
simplicial complex with $I_\Gamma=I(\Delta)$. By using Proposition
\ref{prime} and the previous lemma, we have:

\begin{Corollary}
\label{remember} A subset $F$ of $[n]$ is a facet of $\Gamma$ if
and only if $F^c$ is a minimal vertex cover of $\Delta$.
\end{Corollary}

We say an ideal $I$ in a ring $R$ is {\em Cohen--Macaulay} if
$R/I$ is a Cohen--Macaulay $R$-module. Let $\Delta$ be a
simplicial complex such that $I(\Delta)$ is a Cohen--Macaulay
ideal. Since any Cohen--Macaulay simplicial complex is pure, using
Corollary\ref{remember}, we have $\Delta$ is unmixed.

\medskip

Let $I$ be a squarefree monomial ideal. Then $I=I_\Delta$ for some
simplicial complex $\Delta$. For the convenience we write $I^*$
for $I_{\Delta^\vee}$.

\begin{Lemma}
\label{stra} Let $I$ and $J$ be two squarefree monomial ideals.
Then $$(I\cap J)^*=I^*+J^*.$$
\end{Lemma}
\begin{proof}
Let $P$ be a monomial prime ideal in $R$. Then $I\cap J\subseteq
P$ if and only if $I\subseteq P$ or $J\subseteq P$. The assertion
follows from (\ref{cover}).
\end{proof}

The following theorem gives important algebraic properties of
Alexander duality.

\begin{Theorem}
\label{duality} Let $K$ be a field, $\Delta$ a simplicial complex,
$I_\Delta$ the Stanley--Reisner ideal and $K[\Delta]$ the
Stanley--Reisner ring of $\Delta$. Then
\begin{enumerate}
\item[(1)] {\em (Eagon--Reiner \cite{ER})} $K[\Delta]$ is
Cohen--Macaulay \iff $I_{\Delta^\vee}$ has a linear resolution.

\item[(2)] {\em (Herzog--Hibi--Zheng \cite{HHZ1})} $\Delta$ is
shellable \iff $I_{\Delta^\vee}$ has linear quotients.
\end{enumerate}
\end{Theorem}

\section{A class of unmixed simplicial complexes}

 A simplicial complex $\Delta$ on the vertex set $[n]$ is {\em
Cohen--Macaulay} over a field $K$, if the Stanley--Reisner ideal
$I_\Delta$ of $\Delta$ is a Cohen--Macaulay ideal, while  for a
graph $G$,  we say $G$ is Cohen--Macaulay, if the edge ideal
$I(G)$ of $G$ is a Cohen--Macaulay ideal.

A graph $G$ is {\em bipartite} if its vertex set $V$ can be
partitioned into disjoint subsets $V_1$ and $V_2$ such that every
edge $\{v_1,v_2\}$ of $G$ satisfies $v_1\in V_1$ and $v_2\in V_2$.
Let $G$ be a bipartite graph with no isolated vertex on the vertex
set $V\union V'$, where $V\sect V'=\emptyset$ and $|V|=|V'|$. In
\cite[Theorem 2.4]{HH}, the authors showed that a bipartite graph
$G$ is a Cohen--Macaulay if and only if $I(G)=H^*_\LL$ for some
distributive lattice $\LL$. Later in \cite{HHZ2}, the authors
considered simplicial complexes $\Delta$ on the vertex set
$V\union V'$ with $V\sect V'=\emptyset$ and $|V|=|V'|$, such that
\begin{enumerate}
\item[(1)] there is no $F\in \mathcal{F}(\Delta)$ with $F\subset
V$, and

\item[(2)] $G=\{F\in\mathcal{F}(\Delta)\: F\sect V\neq
\emptyset,\quad F\sect V'\neq \emptyset\}$ is a Cohen--Macaulay
bipartite graph with no isolated vertex,
\end{enumerate}
and showed when the facet ideal $I(\Delta)$ of $\Delta$ is
Cohen--Macaulay, see \cite[Theorem 4.3]{HHZ2}.

In this section we will consider a further  generalization of Theorem
2.4 in \cite{HH}.  For this we need some preparation.

\medskip
The poset ideals and poset coideals of lattices are special
subsets of lattices. Now we introduce a more general class of
subsets of lattices:

\begin{Definition}
\label{segment} {\em Let $\LL$ be a lattice. A subset $\SS$ of
$\LL$ is called a {\em segment} of $\LL$, if  for all $p,q\in\SS$
with $p\leq q$, we have $[p,q]\subseteq \SS$.}
\end{Definition}

It is clear that any poset ideal and any poset coideal of a
lattice $\LL$ are segments of $\LL$. Furthermore, we have

\begin{Lemma}
\label{intersection} Let $\LL$ be a lattice, $\SS$ a subset of
$\LL$. Then the following statements are equivalent:
\begin{enumerate}
\item[(1)] $\SS$ is a segment of $\LL$;

\item[(2)] $\SS$ is the intersection of a poset ideal and a poset
coideal of $\LL$.
\end{enumerate}
\end{Lemma}

\begin{proof}
(1)\implies (2): Let $\II=\{r\in\LL : \text{ there exists an
element }s\in \SS \text { such that } r\leq s\}$ and
$\JJ=\{r\in\LL : \text{ there exists an element } s\in \SS \text {
such that } r\geq s\}.$ Then  $\II$ is a poset ideal of $\LL$ and
$\JJ$ is a poset coideal of $\LL$. For any $s\in \SS$, we have
$s\in \II\cap \JJ$. This implies $\SS\subseteq \II\cap \JJ$. Now
let $r$ be an arbitrary element in $\II\cap \JJ$. Then there exist
$p,q\in \SS$ such that $p\leq r\leq q$, i.e., $r\in[p,q]$. Since
$\SS$ is a segment, we have $r\in \SS$. Hence $\II\cap
\JJ\subseteq \SS$.

(2)\implies (1): Assume $\SS=\II\cap \JJ$, where $\II$ is a poset
ideal of $\LL$ and $\JJ$ is a poset coideal of $\LL$. Let $r\in
[p,q]$ with $p,q\in \SS$ and $p\leq q$. Since $q\in \II$ and
$r\leq q$, we have $r\in \II$. Since $p\in \JJ$ and $r\geq p$, we
have $r\in \JJ$. Hence $r\in \II\sect\JJ=\SS$. This implies that
$\SS$ is a segment of $\LL$.
\end{proof}

\begin{Remark}
\label{minimal} {\em  Let $\SS$ be a segment of a lattice $\LL$.
The poset ideal $\II$ and poset coideal $\JJ$ with the property
$\SS=\II\sect \JJ$ are not uniquely determined.  The poset ideal
$\II$ and poset coideal $\JJ$ in the proof (1)\implies (2) of
Lemma \ref{intersection} are the minimal one with this property.}
\end{Remark}

Let $G$ be a Cohen--Macaulay bipartite graph on the vertex set
$V\union V'$ with $V\sect V'=\emptyset$ and $|V|=|V'|=n$, and
$S=K[x_1,\ldots,x_n,y_1,\ldots,y_n]$ the polynomial ring over a
field $K$. Recall from \cite[Theorem 2.4]{HH} that the vertices
$V=\{x_1,\ldots, x_n\}$ and $V'=\{y_1,\ldots, y_n\}$ can be
labeled such that there exists a partial order $<$ on $V$ with the
property  that $\{x_i,y_j\}$ is an edge of $G$ if and only if
$x_i\leq x_j$. Moreover it is shown that for $P = (V, <)$ the
distributive lattice $\JJ(P)$ satisfies $H_{\JJ(P)}^*=I(G)$. We
denote this lattice by $\LL(G)$. As a  generalization of this
result we have:

\begin{Theorem}
\label{unmixed} Let $\Delta$ be a simplicial complex on the vertex
set $V\union V'$ with $V\sect V'=\emptyset$  and $|V|=|V'|$.
Suppose that $G=\{F\in\mathcal{F}(\Delta)\: F\sect V\neq
\emptyset,\quad F\sect V'\neq \emptyset\}$ is a Cohen--Macaulay
bipartite graph with no isolated vertex. Then the following
conditions are equivalent:
\begin{enumerate}
\item[(1)] $\Delta$ is unmixed, and all minimal vertex covers of
$\Delta$ have cardinality $|V|$;

\item[(2)] there exists a lattice segment $\SS\subseteq \LL(G)$
such that $H_{\SS}^*=I(\Delta)$.
\end{enumerate}
\end{Theorem}

\begin{proof}
(1)\implies (2): Let $\Gamma$ be the (unique) simplicial complex
defined by the equation $I_\Gamma=I(\Delta)$. Since $\Delta$ is
unmixed, we have $\Gamma$ is pure. Let $V=\{x_1,\ldots,x_n\}$ and
$V'=\{y_1,\ldots,y_n\}$ with the labeling as described before this
theorem. Since $\Delta$ is a complex with $2n$ vertices and the
minimal vertex cover of $\Delta$ has cardinality $n$, it follows
from Corollary \ref{remember}, that $|F|=n$ for each
$F\in{\mathcal F}(\Gamma)$.

Let $\Gamma_0$ be the simplicial complex on $V \union V'$ with
$I_{\Gamma_0}=I(G)$. Then  any minimal vertex cover of $\Delta$ is
a minimal vertex cover of $G$. Indeed, a minimal vertex cover $C$
of $\Delta$ is also a vertex cover of $G$, and it has cardinality
$n$, by assumption. On the other hand, since $G$ contains all the
edges $\{x_i,y_i\}$, each vertex cover of $G$ has at least
cardinality $n$. Hence $C$ is a minimal vertex cover  of $G$.

It follows that each facet of $\Gamma$ is a facet of $\Gamma_0$.
In other words, each minimal nonface of $\Gamma^\vee$ is a minimal
nonface of $\Gamma_0^\vee$. Therefore, $G(I_{\Gamma^\vee})\subset
G(I_{\Gamma_0^\vee})=G(H_{\LL(G)})$.  That is, there exists a
subset $\SS\neq\emptyset$ of $\LL(G)$, such that
$G(I_{\Gamma^\vee})=\{u_s :\, s\in \SS\}$, and this implies that
$I(\Delta)=H_\SS^*$.

Now, what we must prove is that for any $p,q\in \SS$ with $p\leq
q$ one has $[p,q]\subseteq \SS$. Suppose, on the contrary, there
exist two elements $\delta$ and $\xi$ of $\LL(G)$ with
$\xi<\delta$, and $\gamma\in \LL(G)$ such that
$\gamma\in[\xi,\delta]$ but $\gamma\notin \SS$.

Recall that the elements  of $\LL(G)$ are poset ideals of
$P=(V,<)$. To simplify the notation, we will assume that $\xi
=\{x_1,\ldots, x_l\}$, $\gamma=\{x_1,\ldots,x_r\}$ and
$\delta=\{x_1,\ldots,x_k\}$ with  $l<r<k$. Since $\xi
=\{x_1,\ldots, x_l\}\in\SS$, we have $x_1\cdots x_ly_{l+1}\cdots
y_n\in G(H_\SS)$. Thus $\{x_1,\ldots,x_l,y_{l+1},\ldots,y_n\}$ is
a minimal vertex cover of $\Delta$. It follows from Corollary
\ref{remember} that $\{y_1, \ldots,y_l,x_{l+1},\ldots,x_n\}\in
{\mathcal F}(\Gamma)$. By the same reason we have
$\{y_1,\ldots,y_k,x_{k+1},\ldots,x_n \}\in {\mathcal F}(\Gamma)$,
but $\{y_1,\ldots,y_r, x_{r+1},\ldots, x_n \}\notin {\mathcal
F}(\Gamma)$. Hence there exists a monomial generator $u$ of
$I_\Gamma=I(\Delta)$ such that $u$ does not divide $y_1\cdots
y_kx_{k+1}\cdots x_n$ and $y_1\cdots y_lx_{l+1}\cdots x_n$,  but
divides $y_1\cdots y_rx_{r+1}\cdots x_n$. Hence there exists an
$i$ with $r<i\leq k$, such that $x_i\mid u$ and a $j$ with
$l<j\leq r$ such that $y_j\mid u$. By our assumption, $u=x_iy_j$.
By our labeling of the vertices it follows that $x_i < x_j$ in
$P$. Since $j\leq r$, we have that $x_j\in \gamma$. Since $\gamma$
is a poset ideal it follows that also $x_i\in \gamma$. This is
impossible, since $i>r$.

(2)\implies (1): Since all generators of $H_\SS$ are of same
degree and since $H^*_\SS=I(\Delta)$, it follows that $\Delta$ is
unmixed.
\end{proof}

The following two simplicial complexes are unmixed and satisfy the
assumption in Theorem \ref{unmixed} with $V_\Delta=\{a,b,c,d\}$,
$V'_\Delta=\{u,v,w,x\}$, and $V_{\Delta'}=\{a,b,c\}$,
$V'_{\Delta'}=\{u,v,w\}$.

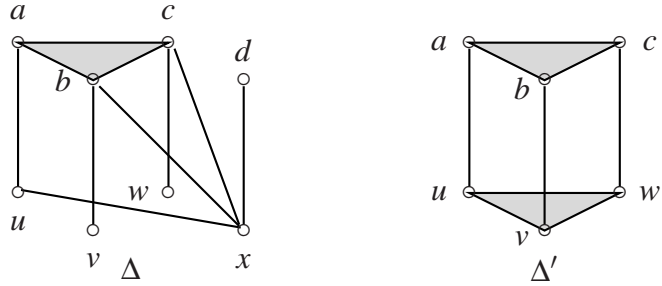
\begin{center}
 \psset{unit=2cm}
\begin{pspicture}(0,-0.3)(5,1.8)
 \psline(0.5,0.54)(0.5,1.46)
 \psline(1,0.29)(1,1.21)
 \psline(1.5,0.54)(1.5,1.46)
 \psline(2,0.29)(2,1.21)
 \psline(0.52,0.52)(1.98,0.27)
 \psline(1.04,1.21)(1.97,0.28)
 \psline(1.54,1.46)(1.97,0.29)
 \pspolygon[style=fyp, fillcolor=medium](0.5,1.5)(1.5,1.5)(1,1.25)
 \rput(0.5,0.5){$\circ$}
 \rput(0.5,1.5){$\circ$}
 \rput(1,0.25){$\circ$}
 \rput(1,1.25){$\circ$}
 \rput(1.5,0.5){$\circ$}
 \rput(1.5,1.5){$\circ$}
 \rput(2,0.25){$\circ$}
 \rput(2,1.25){$\circ$}
 \rput(0.5,0.3){$u$}
 \rput(0.5,1.7){$a$}
 \rput(1,0.05){$v$}
 \rput(0.8,1.25){$b$}
 \rput(1.3,0.5){$w$}
 \rput(1.5,1.7){$c$}
 \rput(2,0.05){$x$}
 \rput(2,1.45){$d$}
 \rput(1.25,0){$\Delta$}

 \psline(3.5,0.54)(3.5,1.46)
 \psline(4.5,0.54)(4.5,1.46)
 \pspolygon[style=fyp, fillcolor=medium](3.5,0.5)(4,0.25)(4.5,0.5)
 \pspolygon[style=fyp, fillcolor=medium](3.5,1.5)(4,1.25)(4.5,1.5)
 \psline(4,0.29)(4,1.21)
 \rput(3.5,0.5){$\circ$}
 \rput(3.5,1.5){$\circ$}
 \rput(4,0.25){$\circ$}
 \rput(4,1.25){$\circ$}
 \rput(4.5,0.5){$\circ$}
 \rput(4.5,1.5){$\circ$}
 \rput(3.3,0.5){$u$}
 \rput(3.3,1.5){$a$}
 \rput(3.85,0.2){$v$}
 \rput(3.85,1.2){$b$}
 \rput(4.7,0.5){$w$}
 \rput(4.7,1.5){$c$}
 \rput(4,0){$\Delta'$}
\end{pspicture}
\end{center}

The segments $\SS$ and $\SS'$ such that $I(\Delta)=H_\SS^*$ and
$I(\Delta')=H_{\SS'}^*$ are given in the next figures. For the
simplicity, sets in these figures are written as monomials. For
example $abcd$ stands for $\{a,b,c,d\}$. The elements of the
segments are indicated by the bullet vertices.

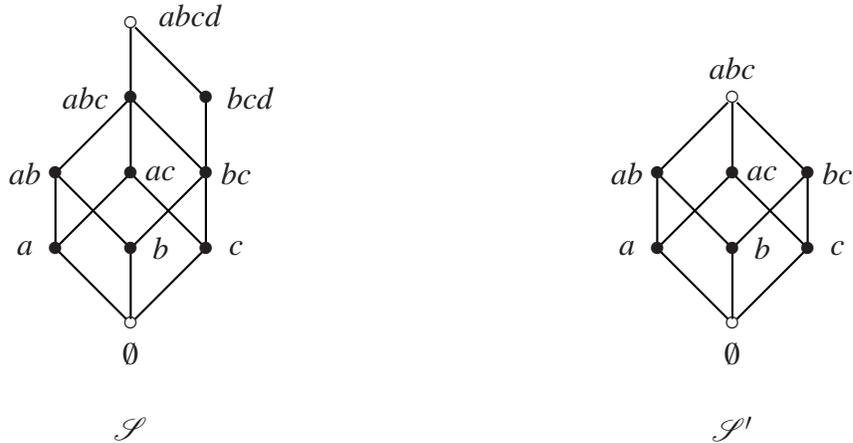
\begin{center}
 \psset{unit=2cm}
\begin{pspicture}(-1,-1)(7,2.3)
\psline(0.97,0.03)(0.5,0.5)
 \psline(0.5,0.5)(1,1)
 \psline(1,1)(1.5,0.5)
 \psline(1.5,0.5)(1.03,0.03)
 \psline(0.5,0.5)(0.5,1)
 \psline(0.5,1)(1,1.5)
 \psline(1,1.5)(1,1)
 \psline(1,1.5)(1.5,1)
 \psline(1.5,1)(1.5,0.5)
 \psline(1,0.03)(1,0.5)
 \psline(0.5,1)(1,0.5)
 \psline(1,0.5)(1.5,1)
 \psline(1,1.5)(1,1.97)
 \psline(1.5,1)(1.5,1.5)
 \psline(1.03,1.97)(1.5,1.5)
 \rput(1,0){$\circ$}
 \rput(0.5,0.5){$\bullet$}
 \rput(1,0.5){$\bullet$}
 \rput(1,2){$\circ$}
 \rput(1.5,1.5){$\bullet$}
 \rput(0.5,1){$\bullet$}
 \rput(1,1){$\bullet$}
 \rput(1.5,0.5){$\bullet$}
 \rput(1.5,1){$\bullet$}
 \rput(1,1.5){$\bullet$}
 \rput(1,-0.2){$\emptyset$}
 \rput(0.3,0.5){$a$}
 \rput(1.2,0.5){$b$}
 \rput(1.7,0.5){$c$}
 \rput(0.3,1){$ab$}
 \rput(1.2,1){$ac$}
 \rput(1.7,1){$bc$}
 \rput(0.7,1.5){$abc$}
 \rput(1.8,1.5){$bcd$}
 \rput(1.4,2.05){$abcd$}
 \rput(1,-0.7){$\SS$}

 \psline(4.97,0.03)(4.5,0.5)
 \psline(4.5,0.5)(5,1)
 \psline(5,1)(5.5,0.5)
 \psline(5,0.03)(5,0.5)
 \psline(4.5,1)(5,0.5)
 \psline(5,0.5)(5.5,1)
 \psline(5.5,0.5)(5.03,0.03)
 \psline(4.5,0.5)(4.5,1)
 \psline(4.5,1)(4.97,1.47)
 \psline(5,1.47)(5,1)
 \psline(5.03,1.47)(5.5,1)
 \psline(5.5,1)(5.5,0.5)
 \rput(5,0){$\circ$}
 \rput(4.5,0.5){$\bullet$}
 \rput(5,0.5){$\bullet$}
 \rput(4.5,1){$\bullet$}
 \rput(5,1){$\bullet$}
 \rput(5.5,0.5){$\bullet$}
 \rput(5.5,1){$\bullet$}
 \rput(5,1.5){$\circ$}
 \rput(5,-0.2){$\emptyset$}
 \rput(4.3,0.5){$a$}
 \rput(5.2,0.5){$b$}
 \rput(5.7,0.5){$c$}
 \rput(4.3,1){$ab$}
 \rput(5.2,1){$ac$}
 \rput(5.7,1){$bc$}
 \rput(5,1.7){$abc$}
 \rput(5,-0.7){$\SS'$}
\end{pspicture}
\end{center}

Note that the facet ideal $I(\Delta)$ of $\delta$ is Cohen--Macaulay,
while the facet ideal $I(\Delta')$ of $\Delta'$ is not. This is
because the ideal $H_\SS=\{avwx,buwx,cuvx,abwx,acvx,bcux,abcx,bcdu\}$ has a linear resolution,
while the ideal $H_{\SS'}=\{avw,buw,cuv,abw,acv,bcu\}$ has no
linear resolution. It is therefore of interest to know for which
kind of segments $\SS$ of a finite distributive lattice $\LL$, the
ideal $H_\SS$ has a linear resolution.

\section{Lattice segments and poset ideals}

We use the notation as in the previous sections.  We have already
seen that $\SS=\II\sect \JJ$ where $\II$ is a poset ideal and
$\JJ$ a poset coideal in $\LL$. In  case $H_\SS=H_\II\sect H_\JJ$,
necessary and sufficient conditions for $H_\SS$ to have a linear
resolution  will be given. We will also discuss when
$H_\SS=H_\II\sect H_\JJ$.

Let $p\in \LL$, and set $N(p)$ for the set of lower neighbors, and
$M(p)$ for the set of upper neighbors of $p$.

Let $P$ be the set of join-irreducible elements of $\LL$, and $<$
a total order on $\LL$ which extends the partial order on $P$. For
a subset $T\subset P$ and $q\in P$ we set
\[
\lambda(q;T)=|\{r\in T\: r<q\}|.
\]
For each element $q\in N(p)$, we have $|\ell(p)\setminus
\ell(q)|=1$. We denote the unique element in $\ell(p)\setminus
\ell(q)$ by $p\setminus q$. Let $p\in \LL$, $S\subseteq N(p)$ and
$T\subseteq M(p)$. We set  $p\setminus S=\{p\setminus s :\, s\in
S\}$,  and $T\setminus p=\{t\setminus p :\, t\in T\}$. Note that
$p\setminus S$ and $T\setminus p$ both are subset of $P$. Let
$p\in \LL$ and $S\subset \LL$. We also set  $p\vee S=\{p\vee s :\,
s\in S\}$ and $p\wedge S=\{p\wedge s :\, s\in S\}$. The following
theorem is shown in \cite{HHZ2}:

\begin{Theorem}[Herzog-Hibi-Zheng \cite{HHZ2}]
 \label{theorem 2.1}
 Let $\LL$ be finite meet-semilattice.
 \begin{enumerate}
 \item[(1)] There exists a finite multigraded free $S$-resolution $\FFF$ of
$H_\LL$ such that for each $i\geq 0$, the free module $\FFF_i$ has
a basis with basis elements
\[
b(p;S)
\]
where $p\in\LL$ and $S\subseteq N(p)$ with $|S|=i$. The
multidegree of $b(p;S)$ is the least common multiple of $u_p$ and
all monomials $u_q$ with $q\in S$.

\item [(2)] The following conditions are equivalent:
\begin{enumerate}
\item [(a)] the resolution constructed in {\em (1)} is minimal;
\item[(b)] for any $p\in\LL$ and for any proper subset $S\subset
N(p)$ the meet $\bigwedge\{q\: q\in S\}$ is strictly greater than
the meet $\bigwedge\{q\:q\in N(p)\}$.
\end{enumerate}

\item[(3)] If $\FFF$ is minimal, then the differential $\partial$
in $\FFF$ is as follows: for each $p\in \LL$ and each $S\subset
N(p)$, one has
\[
\partial(b(p;S))=\sum_{q\in S}(-1)^{\lambda(p\setminus
q; p\setminus S)}(y_{p\setminus
q}b(p;S\setminus\{q\})-x_{p\setminus
q}b(q;q\wedge(S\setminus\{q\})).
\]
\end{enumerate}
 \end{Theorem}

\begin{Corollary}
\label{poset ideal} Let $\LL$ be a finite distributive lattice and
$\II$ a poset ideal (coideal) of $\LL$. Then the minimal free
resolution of $H_\II$ is linear.
\end{Corollary}

\begin{proof}
Let $\II$ be a poset ideal of $\LL$. Then $\II$ is a
meet-semilattice, and has property (2)(b) of Theorem \ref{theorem
2.1}. Hence the free resolution of $H_\II$ as described in Theorem
\ref{theorem 2.1}(1) is minimal. For any $p\in\II$ and any
$S\subseteq N(p)$, the total degree of $b(p;S)$ equals  $\rank \LL
+|S|$. This shows that the resolution of $H_\II$ is linear.

Now assume that $\II$ is a poset coideal. Then by Remark
\ref{veryeasy}, $\widetilde{\II}$ is a poset ideal in
$\widetilde{\LL}$. Therefore $H_{\widetilde{\II}}$  has a linear
resolution by the first part of the proof. By Lemma
\ref{duallattice} (and its proof) the canonical labeling
$\widetilde{\ell}$ of $\widetilde{\LL}$ is given by
$\widetilde{\ell}(p)=P\setminus \ell(p)$ for all
$p\in\widetilde{\LL}$. It follows that $H_{\widetilde{\II}}$ is
generated by the monomials
$\widetilde{u}_p=x_{\widetilde{\ell}(p)}y_{P\setminus
\widetilde{\ell}(p)}=x_{P\setminus \ell(p)}y_{\ell(p)}$. Now we
apply the following involution
\begin{eqnarray}
\label{involution}
 \sigma\: K[\{x_p,y_p\}_{p\in P}]\to K[\{x_p,y_p\}_{p\in
P}], \quad x_p\mapsto y_p\quad \text{and}\quad y_p\mapsto x_p,
\end{eqnarray}
and we obtain  $\sigma(H_{\widetilde{\II}})=H_\II$. This shows
that $H_\II$ has a linear resolution, too.
\end{proof}

As we have already seen that for a poset ideal $\II$ and a poset
coideal $\JJ$ of a finite distributive lattice $\LL$, the ideal
$H_\II$ and $H_\JJ$ both have linear resolutions, one might except
that if we  write $\SS=\II\sect \JJ$ for some poset ideal $\II$
and some poset coideal $\JJ$ of $\LL$ and if $H_\SS=H_\II\sect
H_\JJ$, then the ideal $H_\SS$ has a linear resolution. However
there are two questions arising: (1) when
$H_{\II\sect\JJ}=H_\II\sect H_\JJ$ and (2) whether $H_\II\sect
H_\JJ$ has a linear resolution. In general, the intersection of
two ideals with linear resolutions need not  to have a linear
resolution, even  for the special ideals $H_\II$ and $H_\JJ$. For
example, consider the Boolean lattice $\BB_3$ of rank $3$, and let
$\II=\BB_3\setminus \{\hat 1\}$ and $\JJ=\BB_3\setminus \{\hat
0\}$. Then $H_{\II\sect\JJ}=H_\II\sect H_\JJ$, but it has no
linear resolution.

\medskip
To see when $H_{\II\sect\JJ}=H_\II\sect H_\JJ$ and when
$H_\II\sect H_\JJ$ has a linear resolution, we need some
preparation.

Let $\II$ be a poset ideal of a finite distributive lattice $\LL$
and let $\FFF$ be the minimal free resolution of $H_\LL$ and $\PP$
the minimal free resolution of $H_\II$. Then by Theorem
\ref{theorem 2.1} one sees that $\PP$ is a subcomplex of $\FFF$.
More precisely we have

\begin{Lemma}
\label{basis} For each $i\geq 0$ there exists an injective map
$\Tor_i(K,H_\II)\to \Tor_i(K,H_\LL)$ which maps the basis elements
$b(p;S)$ of $\Tor_i(K,H_\II)$ to the corresponding basis elements
of $\Tor_i(K,H_\LL)$.
\end{Lemma}

\noindent For convenience, in this lemma  and the remaining of
this section the basis elements in a free  resolution and
corresponding basis in the Tor-groups are denoted  by the same
symbol.

\bigskip
Let $\JJ$ be a poset coideal of $\LL$. Then $\widetilde{\JJ}$ is a
poset ideal of $\widetilde{\LL}$. Let $\widetilde{\FFF}$ be the
minimal multigraded free resolution of $H_{\widetilde{\LL}}$ and
$\widetilde{\TT}$ the minimal multigraded free resolution of
$H_{\widetilde{\JJ}}$, as described in Theorem \ref{theorem 2.1}.
Then $\widetilde{\TT}$ is a subcomplex of $\widetilde{\FFF}$, and
the injective map $\Tor_i(K,H_{\widetilde{\JJ}})\to
\Tor_i(K,H_{\widetilde{\LL}})$ is as described in Lemma
\ref{basis}. Since $\sigma(H_{\widetilde{\LL}})=H_\LL$, we have
$\sigma(\widetilde{\FFF})$ is a minimal multigraded free
resolution of $H_\LL$. Since $\FFF$ is also a minimal multigraded
free resolution of $H_\LL$, it is natural to ask what is the
isomorphic chain map from $\widetilde{\FFF}$ to $\FFF$.

To answer this question we need the following two lemmata:

\begin{Lemma}
\label{dual} Let $\widetilde{\LL}$ be the dual of the distributive
lattice $\LL$ and $\widetilde{\FFF}$ the minimal multigraded free
resolution of $H_{\widetilde{\LL}}$ as in Theorem \ref{theorem
2.1}. Then
\begin{enumerate}
\item[(1)] for each $i\geq 0$, the free module $\widetilde F_i$
has a basis with basis elements $\widetilde{b}(r;T)$ with
$r\in\LL$, $T\subseteq M(r)$ in $\LL$ and $|T|=i$. The multidegree
of $\widetilde{b}(r;T)$ is the least common multiple of
$\widetilde{u}_r$ and all monomials $\widetilde{u}_s$ with $s\in
T$;

\item[(2)] the differential $\widetilde{\partial}$ in
$\widetilde{\FFF}$ is as follows: for each $r\in \LL$ and each
$T\subseteq M(r)$, one has
\[
\widetilde{\partial}(\widetilde{b}(r;T))=\sum_{s\in
T}(-1)^{\lambda(s\setminus r;T\setminus r)}(y_{s\setminus
r}\widetilde{b}(r;T\setminus\{s\})-x_{s\setminus
r}\widetilde{b}(s;s\vee(T\setminus\{s\}))).
\]
\end{enumerate}
\end{Lemma}

\begin{proof}
We may assume that $\widetilde{\FFF}$ is a minimal free resolution
of $H_{\widetilde{\LL}}$ as described in Theorem \ref{theorem
2.1}. Therefore $\widetilde{\FFF}$ has a basis
$\widetilde{b}(r;T)$ where $r\in\widetilde{\LL}$ and $T$ is a
subset of lower neighbors of $r$ in $\widetilde{\LL}$. Moreover,
we have
\[
\widetilde{\partial}(\widetilde{b}(r;T))=\sum_{s\in
T}(-1)^{\lambda(\widetilde{r \setminus s};\widetilde{r\setminus
T})}(y_{\widetilde{\ell}(r)\setminus
\widetilde{\ell}(s)}\widetilde{b}(r;T\setminus\{s\})-x_{\widetilde{\ell}(r)\setminus
\widetilde{\ell}(s)}\widetilde{b}(s;s\wedge(T\setminus\{s\})),
\]
where $\widetilde{r \setminus s}$ denote the unique element in
$\widetilde{\ell}(r)\setminus \widetilde{\ell}(s)$ and
$\widetilde{r\setminus T}$ the set $\{\widetilde{\ell}(r)\setminus
\widetilde{\ell}(s) :\, s\in T\}$.

Notice that for any element $r\in \widetilde{\LL}$ (hence $r\in
\LL$, too), a lower (upper) neighbor of $r$ in $\widetilde{\LL}$
is just a upper (lower) neighbor of $r$ in $\LL$, and for any two
element $p$ and $q$ in $\widetilde{\LL}$, the meet (join) of $p$
and $q$ in $\widetilde{\LL}$ is just the join (meet) of them in
$\LL$.

Let $r\in \widetilde{\LL}$ and $s$ a lower neighbor of $r$ in
$\widetilde{\LL}$. We have
\[
\widetilde{\ell}(r)\setminus \widetilde{\ell}(s)=(P\setminus
\ell(r))\setminus (P\setminus \ell(s))=s\setminus r
\]
and
\[
\lambda(\widetilde{r \setminus s};\widetilde{r\setminus
T})=\lambda(\widetilde{\ell}(r)\setminus
\widetilde{\ell}(s);\{\widetilde{\ell}(r)\setminus
\widetilde{\ell}(s) :\, s\in T\})=\lambda(s\setminus r;
\{s\setminus r :\, s\in T\})=\lambda(s\setminus r; T\setminus r).
\]
Thus we obtain the desired formula.
\end{proof}

Let $S$ be any subset of $\LL$. We set $\vee S=\vee\{s:\, s\in
S\}$ and  $\wedge S=\wedge \{s:\, s\in S\}$.

\begin{Lemma}
\label{lcm} Let $\LL$ be a finite distributive lattice, $p\in \LL$
and $S\subseteq N(p)$ with $|S|=i$. Let $r=\wedge\{q:\,q\in S\}$,
and $T$ the set of all upper neighbors of $r$ in the interval
$[r,p]$. Then
\begin{enumerate}
\item[(1)] $|T|=i$ and $\vee T=p$;

\item[(2)] $\lcm(u_p,\{u_q\}_{q\in S})=\lcm(u_r,\{u_s\}_{s\in
T})$;

\item[(3)] for any $r'\neq r$ and $T'\subseteq M(r')$, one has
$\lcm(u_{r'},\{u_{s'}\}_{s'\in T'})\neq \lcm(u_p,\{u_q\}_{q\in
S})$.
\end{enumerate}
\end{Lemma}

\begin{proof}
(1) Since $\LL$ is a distributive lattice, the interval $[r,p]$ is
a Boolean lattice. Hence $|T|=|S|=i$ and $\vee T=p$.

(2) The monomial associated to $p$ is
$u_p=x_{\ell(p)}y_{P\setminus \ell(p)}$, where $P$ is the set of
join irreducible elements of $\LL$. Let $q\in N(p)$. Then
$u_q=x_{\ell(p)\setminus(p\setminus q)}y_{(P\setminus \ell(p))\cup
(p\setminus q)}.$ Hence $$\lcm(u_p,\{u_q\}_{q\in
S})=x_{\ell(p)}y_{(P\setminus \ell(p))\cup(\Union_{q\in
S}(p\setminus q))}.$$ On the other hand, $\ell(r)=\ell(p)\setminus
(\Union_{q\in S}(p\setminus q))$, $u_r=x_{\ell(p)\setminus
(\Union_{q\in S}(p\setminus q))}y_{P\setminus (\ell(p)\setminus
\Union_{q\in S}(p\setminus q))}$.  Since
$\ell(p)=\ell(r)\cup(\Union_{s\in T} \ell(s))$, we have
$$\lcm(u_r,\{u_s\}_{s\in T})=x_{\ell(p)}y_{P\setminus
(\ell(p)\setminus (\Union_{q\in S}(p\setminus q)))}.$$ Since
$\ell(p)\subseteq P$ and $\Union_{p\in S}(p\setminus q)\subseteq
P$, we have $$(P\setminus \ell(p))\cup(\Union_{q\in S}(p\setminus
q))=P\setminus (\ell(p)\setminus \Union_{q\in S}(p\setminus q)).$$
Hence (2) follows.

(3) As in  the proof of (2) we see that the $y$-part of
$\lcm(u_{r'},\{u_{s'}\}_{s'\in T'})$ equals the $y$-part of
$u_{r'}$. Since for any $r'\neq r$, we have $\ell(r')\neq
\ell(r)$. The assertion follows from (2).
\end{proof}

We fix some notation. For each element $r\in \LL$ and $T\subseteq
M(r)$, we write $r^T$ for the join of all elements in $T$, and
$T_r$ the set of all lower neighbors of $r^T$ in the interval
$[r,r^T]$.

The polynomial ring $S$ viewed as a $S$-module via the involution
$\sigma: S\to S$ is denoted by $^\sigma S$. Let $\widetilde{\FFF}$
be the minimal free resolution of the ideal $H_{\widetilde{\LL}}$
with basis elements $\widetilde{b}(r;T)$ as described in Lemma
\ref{dual}. Then $\widetilde{\FFF}\bigotimes_S {^\sigma S}$ with
basis elements $\widetilde{b}(r;T)\bigotimes_S 1$ is a minimal
free resolution of $\sigma(H_{\widetilde{\LL}})=H_\LL$. We denote
the complex $\widetilde{\FFF}\bigotimes_S {^\sigma S}$ by
$\sigma(\widetilde{\FFF})$ and the basis elements
$\widetilde{b}(r;T)\bigotimes_S 1$ by
$\sigma(\widetilde{b}(r;T))$.

\begin{Proposition}
\label{iso} Let $\LL$ be a finite distributive lattice, and let
$\FFF$ and $\widetilde{\FFF}$ be the minimal multigraded free
resolutions of $H_\LL$ and $H_{\widetilde{\LL}}$ , respectively.
Then the map $\pi: \sigma(\widetilde{\FFF})\to \FFF$ with
$\pi(\sigma(\widetilde{b}(r;T)))=(-1)^{|T|}b(r^T;T_r)$ is an
isomorphism of complexes.
\end{Proposition}

\begin{proof}
Let $s\in T$ and $q=\vee \{s'\in T :\, s'\neq s\}$.  Let $|T|=i$.
By Lemma \ref{dual}, we have
\begin{eqnarray}
\label{left}
& &\pi_{i-1}(\widetilde{\partial}_i(\sigma(\widetilde{b}(r;T))))\\
\nonumber &=&\pi_{i-1}(\sum_{s\in T}(-1)^{\lambda(s\setminus
r;T\setminus r)}(y_{s\setminus
r}\widetilde{b}(r;T\setminus\{s\})-x_{s\setminus
r}\widetilde{b}(s;s\vee(T\setminus\{s\}))))\\
\nonumber &=&\sum_{s\in T}(-1)^{\lambda(s\setminus r;T\setminus
r)+(i-1)+1}(y_{s\setminus r}b(s^{s\vee (T\setminus \{s\})};(s\vee
(T\setminus \{s\}))_s)-x_{s\setminus r}b(r^{T\setminus \{s\}};
(T\setminus \{s\})_r))\\
\nonumber &=& \sum_{s\in T}(-1)^{\lambda(s\setminus r;T\setminus
r)+i}(y_{s\setminus r}b(s^{s\vee (T\setminus \{s\})};(s\vee
(T\setminus \{s\}))_s)-x_{s\setminus r}b(r^{T\setminus \{s\}};
(T\setminus \{s\})_r)).
\end{eqnarray}
On the other hand since $\LL$ is a distributive lattice, the
interval $[r,r^T]$ is a Boolean lattice. Hence $q\in T_r$,
$r^T\setminus q=s\setminus r$ and $r^T\setminus T_r=T\setminus r$.
Furthermore, we have
\[
(s\vee (T\setminus \{s\}))_s=T_r\setminus q
\]
and
\[
r^{T\setminus \{s\}}=q, \quad (T\setminus
\{s\})_r=q\wedge(T_r\setminus \{q\}).
\]
These facts together with Theorem \ref{theorem 2.1} yields
\begin{eqnarray}
\label{right} & &\partial_i(\pi_i(\sigma(\widetilde{b}(r;T))))
=\partial_i((-1)^ib(r^T;T_r))\\
\nonumber &=&\sum_{q\in T_r}(-1)^{\lambda(r^T\setminus
q;T_r\setminus r_T)+i}(y_{r^T\setminus q}b(r^T;T_r\setminus
q)-x_{r^T\setminus q}b(q;q\wedge(T_r\setminus \{q\}))\\
\nonumber &=& \sum_{s\in T}(-1)^{\lambda(s\setminus r;T\setminus
r)+i}(y_{s\setminus r}b(s^{s\vee (T\setminus \{s\})};(s\vee
(T\setminus \{s\}))_s)-x_{s\setminus r}b(r^{T\setminus \{s\}};
(T\setminus \{s\})_r)).
\end{eqnarray}
>From (\ref{left}) and (\ref{right}) one sees that $\pi$ is an
isomorphism of complexes.
\end{proof}

Let $\widetilde{\TT}$ and $\widetilde{\FFF}$ be the minimal free
resolutions of $H_{\widetilde{\JJ}}$ and $H_{\widetilde{\LL}}$ as
described in Theorem~\ref{theorem 2.1}, and let $\iota:
\widetilde{\TT}\to \widetilde{\FFF}$ be the injective complex
homomorphism which maps the basis elements $\widetilde{b}(r;T)$ of
$\TT$ to the corresponding basis elements of $\FFF$. Then we have
the following sequence of complex homomorphisms:
\[
\begin{CD}
\widetilde{\TT}@>\iota>>\widetilde{\FFF}@>\sigma>>\sigma(\widetilde{\FFF})@>\pi>>\FFF.
\end{CD}
\]
Let $\psi$ be the map from $\Tor(K,H_{\widetilde{\JJ}})$ to
$\Tor(K,H_\LL)$ induced by  $\pi\circ\sigma\circ\iota$.

As a consequence of the previous proposition, we now have:

\begin{Corollary}
\label{subset} For each $i\geq 0$ the map $\psi_i:
\Tor_i(K,H_{\widetilde{\JJ}})\to \Tor_i(K,H_\LL)$ is injective and
maps the basis elements $\widetilde{b}(r;T)$ of
$\Tor_i(K,H_{\widetilde{\JJ}})$ to the basis elements
$(-1)^{|T|}b(r^T; T_r)$ of $\Tor_i(K,H_\LL)$.
\end{Corollary}

Now we are ready to present one of the main results of this
section:

\begin{Theorem}
\label{equal} Let $\LL$ be a finite distributive lattice, $\II$ a
poset ideal and $\JJ$ a poset coideal of $\LL$ such that
$\II\cup\JJ=\LL$. Then $H_{\II\sect\JJ}=H_\II\sect H_\JJ$ if and
only if for each pair $p,q\in \LL$ with $q\in N(p)$, either
$p\in\II$ or $q\in\JJ$.
\end{Theorem}

\begin{proof}
We may assume that $|\LL|>1$, because otherwise the assertions are
trivial.

Notice that $H_{\II\sect\JJ}\subseteq H_\II\sect H_\JJ$ holds
always, and $H_\II\sect H_\JJ\subseteq H_{\II\sect\JJ}$ if and
only if all generators of $H_\II\sect H_\JJ$ have degree $r$,
where $r$ is the rank of $\LL$.

Consider the long exact Tor-sequence
\[
\begin{CD}
\To \Tor_1(K,H_\II)\dirsum \Tor_1(K,H_\JJ)@>\beta_1
>>\Tor_{1}(K,H_\II+H_\JJ)@>\alpha_1>> \Tor_0(K,H_\II\sect
H_\JJ)\\
\To \Tor_0(K,H_\II)\dirsum \Tor_0(K,H_\JJ)@>>>
\Tor_0(K,H_\II+H_\JJ)@>>> 0
\end{CD}
\]
arising from the short exact sequence
\[
0\To H_\II\sect H_\JJ\To H_\II\dirsum H_\JJ\To  H_\II+H_\JJ\To 0.
\]
Since $\LL=\II\union \JJ$, it follows that $H_\LL=H_\II+H_\JJ$,
and since $H_\LL$ has an $r$-linear resolution we have
$\Tor_i(K,H_\LL)_{i+j}= 0$ if $j\neq r$, and
$\Tor_1(K,H_\LL)_{1+r}\neq 0$. Thus we see that all generators of
$H_\II\sect H_\JJ$ have degree $r$ if and only if $\alpha_1$ is a
zero map, i.e, $\beta_1$ is a surjective map.

By Lemma \ref{basis} we have that the $K$-vector space
$\beta_1(\Tor_1(K, H_\II))$ is spanned by the elements
$b(p;\{q\})$ with $p\in \II$ and $q\in N(p)$, and by Corollary
\ref{subset} we have that the vector space
$\beta_1(\Tor_1(K,H_\JJ))$ is spanned by the elements
$b(q^{\{p\}};\{p\}_q)=b(p;\{q\})$ with $q\in \JJ$ and $p\in M(q)$.
It follows that the image of $\beta_1$ is spanned by the subset
\[
B'=\{b(p;\{q\})\: p\in\II \text{ and } \, q\in N(p), \text{ or }\,
q\in\JJ \text{ and } q\in N(p)\}
\]
of the basis
\[
B=\{b(p;\{q\})\: p\in\LL \text{ and }\, q\in N(p)\}
\]
of $\Tor_1(K,H_\LL)$. Therefore, $\beta_1$ is surjective if and
only if $B'=B$. This implies the assertion.
\end{proof}

It is clear that if $\II\sect\JJ\neq\emptyset$ and
$H_{\II\sect\JJ}\neq H_\II\sect H_\JJ$, then not all  generators
of $H_\II\sect H_\JJ$ have  the same degree. Therefore in this
case,  the ideal $H_\II\sect H_\JJ$ has no linear resolution.
However in case $\II\sect \JJ\neq\emptyset$ and
$H_{\II\sect\JJ}=H_\II\sect H_\JJ$ we have

\begin{Theorem}
\label{linear} Let $\LL$ be a finite distributive lattice, $\II$ a
poset ideal and $\JJ$ a poset coideal of $\LL$ such that
$\II\cup\JJ=\LL$. If $H_{\II\sect\JJ}=H_\II\sect H_\JJ$, then the
following statements are equivalent:
\begin{enumerate}
\item[(1)] $H_\II\sect H_\JJ$ has a linear resolution;

\item[(2)] for each element $p\in\LL$, either $p\in\II$ or $\wedge
N(p)\in \JJ$;

\item[(3)] for each element $r\in\LL$, either $r\in\JJ$ or $\vee
M(r)\in \II$.

\end{enumerate}
\end{Theorem}

\begin{proof}
We may assume that $|\LL|>1$.

(2)\implies(3): Assume there exists some element $r\in \LL$ such
that $r\not\in\JJ$ and $p=\vee M(r)$ does not belong to $\II$.
Since $\LL$ is a distributive lattice, the interval $[r,p]$ is a
Boolean lattice. Hence $\wedge N(p)=r$. Therefore we have
$p\not\in \II$ and $r=\wedge N(p)\not\in \JJ$, a contradiction.

By the same argument, one sees that (3) implies (2).

Now, we prove that the conditions (1) and (2) are equivalent.
Consider the long exact Tor-sequence
\[
\begin{CD}
\cdots\to \Tor_{i+1}(K,H_\II)\dirsum
\Tor_{i+1}(K,H_\JJ)@>\beta_{i+1}
>>\Tor_{i+1}(K,H_\LL)@>\alpha_{i+1}>>
\Tor_i(K,H_\II\sect H_\JJ) \\
\To \Tor_i(K,H_\II)\dirsum \Tor_i(K,H_\JJ)@>>> \Tor_i(K,H_\LL)@>>>
\cdots
\end{CD}
\] arising from the short exact sequence
\[
0\To H_\II\sect H_\JJ\To H_\II\dirsum H_\JJ\To  H_\LL\To 0.
\]
Here we used that $\II\union\JJ=\LL$, so that $H_\II+H_\JJ=H_\LL$.
Let $r=\rank \LL$. Since the ideal $H_\II\sect
H_\JJ=H_{\II\sect\JJ}$ is generated in degree $r$, it has a linear
resolution if and only if
\[
\begin{CD}
\Tor_{i+1}(K,H_\LL)_j@>\alpha_{i+1,j}>> \Tor_i(K,H_\II\sect
H_\JJ)_j
\end{CD}
\]
is the zero map for all $j\neq i+r$ and all $i\geq 0$, since the
ideals $H_\II$ and $H_\JJ$ have $r$-linear resolutions. Since the
ideal $H_\LL$ has an $r$-linear resolution, the map
$\alpha_{i+1,j}=0$ for $j\neq i+1+r$. Hence $H_\II\sect H_\JJ$ has
a linear resolution if and only if $\alpha_{i+1,i+1+r}=0$ for all
$i\geq 0$, and this is the case if and only if $\beta_{i+1,i+1+r}$
is surjective for all $i\geq 0$.

We argue as in the proof of Theorem \ref{equal}. The set
\[
B_{i+1}=\{b(p;S)\: p\in\LL \text{ and }\, S\subset N(p),\,
|S|=i+1\}
\]
is a $K$-basis of $\Tor_{i+1}(K,H_\LL)_{i+1+r}$. Using Lemma
\ref{basis} and Corollary \ref{subset} we see that the set
\begin{eqnarray*}
B_{i+1}'&=&\{b(p;S)\: p\in \II  \text{ and }\, S\subset N(p),\,
|S|=i+1\}\\
&\union& \{b(r^T;T_r)\: r\in\JJ \text{ and }\, T\subset M(r),\,
|T|=i+1\} \end{eqnarray*} spans the image of $\beta_{i+1,i+1+r}$.
Thus $\beta_{i+1,i+1+r}$ is surjective if and only if
$B_{i+1}'=B_{i+1}$ for all $i> 0$.  Note that $B_{i+1}'\subset
B_{i+1}$. Suppose condition (2) holds, and let $b(p;S)\in
B_{i+1}$. If $p\in \II$, then $b(p;S)\in B_{i+1}'$. If $p\not\in
\II$, then $p\in \JJ$. Let $r=\wedge S$. It follows from condition
(2) that  $r\in \JJ$. Let $T$ be the set of upper neighbors of $r$
in the interval $[r,p]$. Then $p=r^T$ and $S=T_r$, and hence
$b(p;S)=b(r^T;T_r)$ belongs to $B'_{i+1}$.

Conversely assume that $B_{i+1}'=B_{i+1}$. In particular, for all
$p\in\LL$ we have $b(p;N(p))\in B_{|N(p)|}'$. So either $p\in\II$
or there is some $r\in \JJ$ and $T\subset M(r)$ such that $p=r^T$
and $T_r=N(p)$. Since $\wedge T_r=r$, it follows that $\wedge
N(p)=r$ which is in $\JJ$.
\end{proof}

Up to now, we always assume that $\II\cup \JJ=\LL$ and $\II\sect
\JJ\neq\emptyset$.  Now we consider the case $\II\cup \JJ=\LL$ and
$\II\sect \JJ=\emptyset$. In this case, $H_{\II\sect\JJ}\neq
H_\II\sect H_\JJ$. However, we have:

\begin{Theorem}
\label{empty} Let $\LL$ be a finite distributive lattice, $\II$ a
poset ideal and $\JJ$ a poset coideal of $\LL$. If $\II\cup
\JJ=\LL$ and $\II\sect\JJ=\emptyset$, then the ideal $H_\II\sect
H_\JJ$ has a linear resolution.
\end{Theorem}

This theorem follows immediately from the following two lemmata.

\begin{Lemma}
\label{generator} Under the assumption of Proposition \ref{empty},
the ideal $H_\II\sect H_\JJ$ is generated by the monomials
$\lcm(u_p,u_q)$ with $p\in\JJ$, $q\in \II$ and $q\in N(p)$. In
particular, all generators of $H_\II\sect H_\JJ$ are of degree
$\rank \LL+1$.
\end{Lemma}

\begin{proof}
Since $\II\sect \JJ=\emptyset$, all generators of  $H_\II\sect
H_\JJ$ have degree greater than the rank of $\LL$. Let
$H=\langle\lcm(u_p,u_q):\, p\in\JJ, q\in \II \text{ and } q\in
N(p)\rangle$. It is clear that $H\subseteq H_\II\sect H_\JJ$.
Since $H_\II$ and $H_\JJ$ both are monomial ideals, the
intersection $H_\II\sect H_\JJ$ is again a monomial ideal. Let $m$
be any monomial in $H_\II\sect H_\JJ$. Then there exist $r\in \II$
and $s\in\JJ$, such that $u_r| m$ and $u_s|m$. Let $C$ be any
chain between $r\wedge s$ and $r\vee s$. Since $\II$ is a poset
ideal and $\JJ$ is a poset coideal, and $r\in \II$, $s\in\JJ$, we
have $r\wedge s\in\II$ and $r\vee s\in\JJ$. Hence there exist
$p,q\in C$ such that $q\in\II$, $p\in \JJ$ and $q$ is a lower
neighbor of $p$. We claim $\lcm(u_p,u_q)| m$. To see this, we
write $m=\prod_{i=1}^nx^{a_i}y^{b_i}$ as $m_xm_y$ where
$m_x=\prod_{i=1}^nx^{a_i}$ and $m_y=\prod_{i=1}^ny^{b_i}$,  and as
before we write $u_t=x_{\ell(t)}y_{P\setminus \ell(t)}$, where $P$
is the set of join irreducible elements of $\LL$. Since
$\ell(p)\subseteq\ell(r\vee s)$, we have $x_{\ell(p)}|
x_{\ell(r\vee s)}$. Since $u_r| m$ and $u_s| m$, we have
$x_{\ell(r\vee s)}|m$ and hence  $x_{\ell(p)}| m_x$. By the same
argument we see that $y_{P\setminus {\ell(q)}}| m_y$. Hence
$\lcm(u_p,u_q)=x_{\ell(p)}y_{P\setminus  \ell(q)}$ divides $m$.
\end{proof}

\begin{Lemma}
\label{d+1} Let $R$ be a polynomial ring over a field $K$, $I$ and
$J$ ideals in $R$. Suppose $I$, $J$ and $I+J$ have $d$-linear
resolutions. If all elements of $G(I\sect J)$ have degree $d+1$,
then the ideal $I\sect J$ has a $(d+1)$-linear resolution.
\end{Lemma}

\begin{proof}
Consider the long exact Tor-sequence
\[
\begin{CD}
\cdots\to \Tor_{i+1}(K,I)\dirsum \Tor_{i+1}(K,J)@>\beta_{i+1}
>>\Tor_{i+1}(K,I+J)@>\alpha_{i+1}>>
\Tor_i(K,I\sect J) \\
\To \Tor_i(K,I)\dirsum \Tor_i(K,J)@>>> \Tor_i(K,I+J)@>>>\cdots.
\end{CD}
\]
Since $I$, $J$ and $I+J$ have $d$-linear resolutions. It follows
that $\Tor_i(K,I)_j=\Tor_i(K,J)_j=0$ for any $j\neq i+d$ and
$\Tor_{i+1}(K,I+J)_j=0$ for any $j\neq i+1+d$. Hence
$\Tor_i(K,I\sect J)_j=0$ for $j<i+d$ or $j>i+1+d$. Since $I\sect
J$ is generated in degree $d+1$, we have $\Tor_i(K,I\sect J)_j=0$
for $j=i+d$. Therefore $\Tor_i(K,I\sect J)_j=0$ for any $j\neq
i+1+d$. Hence $I\sect J$ has a $(d+1)$-linear resolution.
\end{proof}

\medskip\noindent
In the remaining of this section we discuss some special classes
of segments of a finite distributive lattice $\LL$ to which our
results apply and where some additional information can be
obtained . As we have already seen, for any segment $\SS$ there
exist a poset ideal $\II$ and a poset coideal $\JJ$ such that
$\SS=\II\sect \JJ$. Now Let $\LL$ be a finite distributive lattice
of rank $r$. We consider a special class of segments $\SS$ of
$\LL$ which consisting of all elements $p$ in $\LL$ such that
$i\leq \rank p\leq j$ for some $i$ and $j$ with $0\leq i\leq j\leq
r$. We denote it by $\LL_{i,j}$.

\begin{Lemma}
\label{choose} Let $\LL$ be a finite distributive lattice of rank
$r$, and let $\LL_{i,j}$ be a segment of $\LL$. Then  there exists
a poset ideal $\II$ and a poset coideal $\JJ$ of $\LL$ such that
$H_{\LL_{i,j}}=H_{\II\sect\JJ}=H_\II\sect H_\JJ$.
\end{Lemma}

\begin{proof}
Let $\II=\{p\in\LL:\, \rank p\leq j\}$, and $\JJ=\{p\in\LL:\,
\rank p\geq i\}$. Then $\II$ is a poset ideal, $\JJ$ is a poset
coideal of $\LL$, and $\LL_{i,j}=\II\sect\JJ$. It remains to show
that $H_{\II\sect\JJ}=H_\II\sect H_\JJ$. Let $p,q\in\LL$ and $q\in
N(p)$. If $p\notin \II$, then $\rank p>j$. Hence $\rank q=\rank
p-1\geq j\geq i$, i.e., $q\in \JJ$. The assertion follows from
Proposition \ref{equal}.
\end{proof}

With the assumptions and notation of the previous lemma, the ideal
$H_{\LL_{i,j}}$ has a linear resolution if and only if $H_\II\sect
H_\JJ$ has a linear resolution.

\begin{Corollary}
\label{top and bottom} Let $\LL\neq\{\hat 0,\hat 1\}$ be a finite
distributive lattice and $\SS=\LL\setminus \{\hat 0,\hat 1\}$.
Then
\begin{enumerate}
\item[(1)] $H_\SS$ has a linear resolution if and only if $\LL$ is
not a Boolean lattice;

\item[(2)] in case $\LL$ is the Boolean lattice $\BB_r$, the ideal
$H_\SS$ has the following minimal free resolution:
\[
\begin{CD}
\TT \: 0\To T_{r-1}\To \cdots T_1\To T_0\To  H_\SS\To 0.
\end{CD}
\]
\end{enumerate}
with $T_i=S^{\binom{r}{i}(2^{r-i}-2)}(-r-i)$ for $i=0,\ldots, r-2$
and $T_{r-1}=S(-2r)$.

\end{Corollary}

\begin{proof}
(1) Let $\II=\{p\in\LL:\, \rank p\leq \rank \LL-1\}$, and
$\JJ=\{p\in\LL:\, \rank p\geq 1\}$. Then by Lemma \ref{choose},
$\II$ and $\JJ$ are the poset ideal and poset coideal of $\LL$
such that $H_\SS=H_\II\sect H_\JJ$. The distributive lattice $\LL$
is a Boolean lattice if and only if the meet of all lower
neighbors of $\hat 1$ is $\hat 0$. Since $\hat 1$ is the only
element which is not in $\II$, and $\hat 0$ is the only element
which is not in $\JJ$, by using the Theorem \ref{linear}, we have
$H_\SS$ has a linear resolution if and only if $\LL$ is not a
Boolean lattice.

(2) Choose $\II$ and $\JJ$ as in the proof of (1). Hence
$H_\SS=H_{\II\sect\JJ}=H_\II\sect H_\JJ$. We consider long exact
Tor-sequence as in the proof of Theorem \ref{linear}. Notice that
$(\hat 1, N(\hat 1))$ is the only pair with the form $(p,S)$ with
$p\in\BB_r$ and $S\subset N(p)$ such that $p\not\in \II$ and
$\wedge S\not\in\JJ$. It follows that  the map $\beta_i$ is
surjective for $i<r$.  Hence  for all $i<r-1$ we have the exact
sequence
\begin{eqnarray}
\label{sooninjapan}
\begin{CD}
\hspace{0.5cm} 0\To\Tor_i(K,H_\SS) \To \Tor_i(K,H_\II)\dirsum
\Tor_i(K,H_\JJ)\To \Tor_{i}(K,H_{\BB_r})\To 0,
\end{CD}
\end{eqnarray}
and so
\[
b_i(H_\SS)=b_i(H_\II)+b_i(H_\JJ)-b_i(H_{\BB_r})\quad
\text{for}\quad i<r-1,
\]
where $b_i(I)$ is the $i$-th Betti number of the ideal $I$. By
using Theorem \ref{theorem 2.1}, Corollary \ref{subset} and the
combinatorial fact that each Boolean lattice $\BB_r$ contains
$\binom{r}{i}2^{r-i}$ Boolean sublattices $\BB_i$, we have
$b_i(H_\II)=b_i(H_\JJ)=\binom{r}{i}2^{r-i}-\binom{r}{i}$ and
$b_i(H_\LL)=\binom{r}{i}2^{r-i}$. Hence
\[
b_i(H_\SS)=2(\binom{r}{i}2^{r-i}-\binom{r}{i})-\binom{r}{i}2^{r-i}=\binom{r}{i}(2^{r-i}-2),
\]
for any $i<r-1$. It also follows from (\ref{sooninjapan}) that the
resolution of $H_\SS$ is linear up to homological degree $r-2$.

Now let $i=r-1$. Since $\Tor_r(K,H_\II)= \Tor_r(K,H_\JJ)=0$,  we
get the exact sequence
\begin{eqnarray*}
0\to \Tor_r(K,H_{\BB_r})\to&\Tor_{r-1}(K,H_\SS) &\to
\Tor_{r-1}(K,H_\II)\dirsum \Tor_{r-1}(K,H_\JJ)\\
&\To &\Tor_{r-1}(K,H_{\BB_r})\to 0.
\end{eqnarray*}
Since $\dim_K \Tor_{r-1}(K,H_\II)\dirsum
\Tor_{r-1}(K,H_\JJ)=\dim_K \Tor_{r-1}(K,H_\JJ)=2r$, it follows
that $\Tor_{r-1}(K,H_\II)\dirsum \Tor_{r-1}(K,H_\JJ) \to
\Tor_{r-1}(K,H_{\BB_r})$ is an isomorphism. Hence
$$\Tor_{r-1}(K,H_\SS)\iso \Tor_r(K,H_{\BB_r})\iso K(-2r),$$
as desired.
\end{proof}

Using Lemma \ref{choose} and Theorem\ref{linear}, we have the
following two facts:

\begin{Corollary}
\label{one line} Let $\LL$ be a finite distributive lattice, and
$i$ an integer. If $|\LL_{i,i}|>1$, then the ideal $H_{\LL_{i,i}}$
has no linear resolution.
\end{Corollary}

\begin{proof}
Let $\II=\{p\in\LL:\, \rank p\leq i\}$ and $\JJ=\{p\in\LL:\, \rank
p\geq i\}$. Then by Lemma \ref{choose} we have
$H_{\LL_{i,i}}=H_{\II\sect\JJ}=H_\II\sect H_\JJ$. Since $\LL$ is
distributive, there exist elements $u$ and $v$ in $\LL_{i,i}$ such
that $\rank (u\vee v)=i+1$. Hence $\rank (\wedge N(u\vee v))<\rank
u=i$. Therefore we have $u\vee v\not\in \II$ and $\wedge N(u\vee
v)\not\in \JJ$. By Theorem \ref{linear}, $H_{\LL_{i,i}}$ has no
linear resolution.
\end{proof}

\begin{Corollary}
\label{planar} If $\LL$ is a finite planar distributive lattice of
rank $r$, then the ideal $H_{\LL_{i,j}}$ has a linear resolution,
if  $i<j$.
\end{Corollary}

\begin{proof}
Let $\II=\{p\in\LL:\, \rank p\leq j\}$ and $\JJ=\{p\in\LL:\, \rank
p\geq i\}$. Then $H_{\LL_{i,j}}=H_{\II\sect\JJ}=H_\II\sect H_\JJ$.
Since $\LL$ is a planar distributive lattice, each element $p$ in
$\LL$ has at most two lower neighbors. Hence $\rank \wedge
N(p)\geq p-2$. Thus if $p\not\in \II$, then $\rank p>j$. Therefore
since $i<j$, we have $\rank q\geq \rank p-2\geq i$, i.e.,
$q\in\JJ$. By Theorem \ref{linear}, $H_{\LL_{i,j}}$ has a linear
resolution.
\end{proof}

\end{document}